\documentclass[11pt]{amsart}

\usepackage{etex}
\usepackage{savesym}
\usepackage{amssymb, amsmath,amsthm,amscd}
\usepackage{dsfont}
\savesymbol{lrcorner}\savesymbol{llcorner}\savesymbol{urcorner}\savesymbol{ulcorner}
\usepackage{mathabx}
\usepackage{stmaryrd}
\usepackage{setspace}
\usepackage{chngcntr}
\usepackage{enumerate,makecell}
\usepackage[usenames,dvipsnames]{xcolor}
\usepackage[bookmarks=true,colorlinks=true, linkcolor={blue!35!black}, citecolor=cyan]{hyperref}
\usepackage{lmodern}
\usepackage{graphicx,float}
\usepackage{multirow}
\usepackage{geometry}
\newgeometry{left=1.45in,right=1.45in,top=1.5in,bottom=1.5in}
\restoresymbol{txfonts}{lrcorner}\restoresymbol{txfonts}{llcorner}\restoresymbol{txfonts}{urcorner}\restoresymbol{txfonts}{ulcorner}

\usepackage{color}
\usepackage[all,poly,color,matrix,arrow]{xy}

\newcommand{\comment}[1]{}

\def\tn{\textnormal}

\def\mc{\mathcal}

\def\RR{{\mathbb R}}

\def\NN{{\mathbb N}}

\def\bigob{\mathds{O}}
\def\bigmor{\mathds{M}}

\def\Ob{\tn{Ob}}

\def\SEL*{\tn{SEL*}}

\def\hsp{\hspace{.3in}}

\def\to{\rightarrow}

\def\taking{\colon}

\def\too{\longrightarrow}

\def\|{{\;|\;}}
\def\m1{{-1}}

\def\ullimit{\ar@{}[rd]|(.3)*+{\lrcorner}}
\def\urlimit{\ar@{}[ld]|(.3)*+{\llcorner}}
\def\lllimit{\ar@{}[ru]|(.3)*+{\urcorner}}
\def\lrlimit{\ar@{}[lu]|(.3)*+{\ulcorner}}
\def\ulhlimit{\ar@{}[rd]|(.3)*+{\diamond}}
\def\urhlimit{\ar@{}[ld]|(.3)*+{\diamond}}
\def\llhlimit{\ar@{}[ru]|(.3)*+{\diamond}}
\def\lrhlimit{\ar@{}[lu]|(.3)*+{\diamond}}
\newcommand{\clabel}[1]{\ar@{}[rd]|(.5)*+{#1}}
\newcommand{\TriRight}[7]{\xymatrix{#1\ar[dr]_{#2}\ar[rr]^{#3}&&#4\ar[dl]^{#5}\\&#6\ar@{}[u] |{\Longrightarrow}\ar@{}[u]|>>>>{#7}}}
\newcommand{\TriLeft}[7]{\xymatrix{#1\ar[dr]_{#2}\ar[rr]^{#3}&&#4\ar[dl]^{#5}\\&#6\ar@{}[u] |{\Longleftarrow}\ar@{}[u]|>>>>{#7}}}
\newcommand{\TriIso}[7]{\xymatrix{#1\ar[dr]_{#2}\ar[rr]^{#3}&&#4\ar[dl]^{#5}\\&#6\ar@{}[u] |{\Longleftrightarrow}\ar@{}[u]|>>>>{#7}}}

\newcommand{\arr}[1]{\ar@<.5ex>[#1]\ar@<-.5ex>[#1]}
\newcommand{\arrr}[1]{\ar@<.7ex>[#1]\ar@<0ex>[#1]\ar@<-.7ex>[#1]}
\newcommand{\arrrr}[1]{\ar@<.8ex>[#1]\ar@<.3ex>[#1]\ar@<-.3ex>[#1]\ar@<-.8ex>[#1]}
\newcommand{\arrrrr}[1]{\ar@<1ex>[#1]\ar@<.5ex>[#1]\ar[#1]\ar@<-.5ex>[#1]\ar@<-1ex>[#1]}

\newcommand{\To}[1]{\xrightarrow{#1}}

\newcommand{\Adjoint}[4]{\xymatrix@1{#2 \ar@<.5ex>[r]^-{#1} & #3 \ar@<.5ex>[l]^-{#4}}}
\newcommand{\adjoint}[4]{\xymatrix{#1\taking #2\ar@<.5ex>[r]& #3\hspace{1pt}:\hspace{-2pt} #4\ar@<.5ex>[l]}}

\def\id{\tn{id}}

\def\Top{\mathsf{Top}}
\def\Cat{\mathsf{Cat}}

\def\Mon{\mathsf{Mon}}
\def\Grp{\mathsf{Grp}}

\def\Vect{\mathsf{Vect}}

\def\Set{\mathsf{Set}}

\def\mcC{\mc{C}}
\def\mcD{\mc{D}}

\def\mcG{\mc{G}}

\theoremstyle{remark}

\newtheorem{question}[subsubsection]{Question}
\newtheorem{guess}[subsubsection]{Guess}

\theoremstyle{definition}

\usepackage{tikz}
\usetikzlibrary{arrows,calc,chains,matrix,positioning,scopes,snakes}

\DeclareMathOperator{\MatSet}{\mathsf{Mat}}
\newcommand{\Mat}[2]{\MatSet_{#1\times #2}}
\newcommand{\MG}[1]{\mathsf{InvMat}_{#1}}
\newcommand{\MS}[1]{\Mat{#1}{#1}}

\begin{document}

\title{Categories as mathematical models}

\author{David I. Spivak}
\address{Massachusetts Institute of Technology, Department of Mathematics, 77 Massachusetts Avenue, Cambridge, MA 02139-4307, USA}
\email{dspivak@gmail.com}

\maketitle

\begin{abstract}
Category theory is presented as a mathematical modeling framework that highlights the relationships between objects, rather than the objects in themselves. A working definition of model is given, and several examples of mathematical objects, such as vector spaces, groups, and dynamical systems, are considered as categorical models.
\end{abstract}

\tableofcontents
\small
\begin{quote}
Mathematicians do not study objects, but relations between objects. Thus, they are free to replace some objects by others so long as the relations remain unchanged. Content to them is irrelevant: they are interested in form only. 
\begin{flushright} -- \emph{Henri Poincar\'e}\end{flushright}
\end{quote}
\bigskip
\begin{quote}
Yet, I hope that I managed to convey the message: the mathematical language developed by the end of the 20th century by far exceeds in its expressive power anything, even imaginable, say, before 1960. Any \emph{meaningful} idea coming from science can be fully developed in this language. 
\begin{flushright} -- \emph{Mikhail Gromov}\end{flushright}
\end{quote}
\normalsize

\section{Introduction}\label{sec:intro}

In the sciences, most of the prominent methods for incorporating mathematics involve setting up stochastic processes, dynamical systems, or statistical models that capture the relevant matters in the scientific subject. At bottom, these techniques all involve interplays of numbers. However, a biologist's sophisticated understanding of the complex aspects of life---heredity, reproduction, hierarchical nesting, symbiosis, metabolism, etc.---remains trapped in the realm of  \emph{ideas}.%
\footnote{Here, we use biology only as a specific branch of science. Our intention is to consider \emph{scientific ideas}, beyond their numerical shadows, as subject to mathematical formalization.}  
 Such ideas can often be reduced to numerical models, but the ideas themselves are confined to the background. This is suboptimal, because when an idea or theory is itself mathematically formalized, it gains clarity, systematicity, and falsifiability. 
%Surely, a biologist's understanding motivates her hypotheses and experimental designs, and these can be refined in her discussions with other biologists. However, unless the ideas themselves can be captured mathematically, they are not part of science for Kant. 

My hope for category theory is that it can be used to model many of the actual ideas and ways of thinking that exist within science. Modeling a phenomenon allows us to examine, and interact with, a simplified version of it, and involving mathematics generally provides an additional level of rigor and communicability. A category-theoretic model of the ideas, rather than of the mere quantities, may be able to formally capture a whole conceptual framework, say, ideas about the hierarchical nature of organizational systems. With such a formal description, one could apply rigorous conceptual---rather than numerical---tests to the ideas themselves. For example, one might check whether the ideas satisfy various internal consistencies, or one might ask about the nature of their integration with the other major conceptual frameworks that exist in the field. 

There is a good deal of work on using category theory to model high-level conceptual aspects of scientific subjects. For example, categories have been used by John Baez to model signal flow and reaction networks, by Abramsky and Coecke to model aspects of quantum mechanics, and by Lambek to model computer programming languages. I believe that in 2015 we are at the early stages of an effort to ``categorify" science.

Category theory holds promise for such a conceptual integration of science because it has achieved such an integration of mathematics to a remarkable degree. It does so, not by finding a single syntax or format that can encode any structure, but by doing the opposite: leaving the encoding entirely absent. Rather than modeling a given object \emph{in itself}, category theory models only the relationships between objects. For example, one may say that some object has an internal structure; however, in the categorical model this internal structure can only be seen by using other objects as probes. That is, category theory views an object's internal structure only in terms of its relationships with other objects.

%As mentioned above, there exist a plethora of mathematical techniques, most of which are numerical in nature,  for modeling scientific subjects. And quite a few of these, including calculus, dynamical systems, logic, and probability theory, have well-established category-theoretic descriptions. However, there are many other categories that are not numerical in nature, such as groups, monoids, and operads, which may also be quite useful in scientific modeling. 

My goal in this paper is to give philosophers an intuitive idea about how category theory can be thought of as a universal modeling language, in which the relationships between objects are paramount. In particular, it is not my goal to present new interpretations of established philosophical concepts. At times I invoke philosophical ideas, e.g., those of Kant, but the reader should keep in mind that these are my attempt to forge a communication with the philosophical community, rather than to express substantial claims about Kant's work. It is also not my goal to explain exactly how category theory can be used to formalize the more complex ideas found throughout science. As mentioned above, this kind of work has certainly begun in earnest, but adequately addressing it is beyond the scope of this paper. The discussion in this paper is somewhat similar in spirit to that of Lawvere and Schanuel's excellent introductory book, \emph{Conceptual Mathematics}, which I recommend for further reading. 

\subsection{Category theory as modeling language}\label{subsec:CT_modeling}

In this paper I cast category theory (CT) as a universal modeling language. More precisely, I claim that CT is \emph{a mathematical model of mathematical models}. Before I explain this assertion, I will ground the discussion with a prototypical example of a model to fix ideas, and then I will define what I mean by model. 

The sonar system in a submarine models the distances between various objects in the physical environment of the submarine. It does so by plotting dots, which correspond to these objects, on a screen using polar coordinates (distance and angle). A representation of the information gathered by the sonar system is displayed to the submarine's pilot in a familiar ``self-centric" way, i.e., with the submarine (and hence the pilot) shown at the center of the display. This allows the sonar system to become \emph{transparent equipment}%
\footnote{See Andy Clark, \emph{Supersizing the Mind}.} 
for the pilot. In other words, she can seamlessly integrate the representation into her personal repertoire and thereby use familiar methods to cope appropriately with a situation as it unfolds, say to evade or pursue another object. 

In this example, the sonar system of the submarine is the \emph{model} and the pilot is its \emph{user}. The physical environment of the submarine is the \emph{subject} of the model; it is the thing being modeled by the sonar system. Because the locations of objects in the physical environment are emphasized, rather than suppressed, we say that these observable aspects of the subject are \emph{foregrounded} by the model. The distances and angles between these locations are the \emph{relationships} between foregrounded aspects. We refer to the translation, from sound-wave reflectance times to points in a polar coordinate system, as the \emph{formalism} that founds the model. 

Let us tie these ideas together by making some observations. The value of the sonar system (the model) is measured by the extent to which the pilot's (the user's) interaction with submarine's physical environment (the subject) is successfully mediated by the sonar system. This in turn depends on the extent to which the locations of the physical objects in the environment (the foregrounded aspects) are propitious for successful interaction with the subject and the extent to which the polar coordinate representation of these objects (the formalism) is faithful.

With this exemplar in mind, we make the following two philosophical postulates, which will help to organize the ideas in this paper. 
\begin{enumerate}[\hspace{.19in} (1)]
\item Modeling a subject is foregrounding certain observable aspects of the subject, and then formalizing these aspects and certain observable relationships between them. 
\item The value of a model is measured by the extent to which the user's interactions with the subject are successfully mediated by the model. This in turn depends on the propitiousness of the foregrounded aspects and the faithfulness of the formalism.
\end{enumerate}
Note that mathematical models seem to put more emphasis on, and care into, the formalism than do other types of models.

\subsection{Using models is connecting models}

In his \emph{Critique of Pure Reason}, Immanuel Kant makes an important assertion:
\begin{quote}[E]verything intuited or perceived in space and time, and therefore all objects of a possible experience, are nothing but phenomenal appearances, that is, mere representations, which in the way in which they are represented to us, as extended beings, or as series of changes, have no independent, self-subsistent existence apart from our thoughts.
\end{quote}

In other words, our interactions with the subject, which is ostensibly {\em out there}, are actually interactions with our own familiar models (or as Kant says, representations) of it. Our mind is an economy of models and our thinking consists of negotiations within that economy. Thus the value of a model is measured by the ease with which it negotiates or interfaces with the other models in our repertoire.

Using models is all about translating between models. To say it another way, the observable aspects of a model are known only by its relationships with other models. It follows that, in order to objectively understand our interactions with a model, it is useful to understand the more general question of how models relate to other models. We add to our list a third and final philosophical postulate, which will be clarified throughout the paper:
\begin{enumerate}[\hspace{.19in} (3)]
\item A model is known only by its relationships with other models.
\end{enumerate}

With our three modeling postulates in hand, we unpack the italicized statement from Section~\ref{subsec:CT_modeling}, that category theory is a mathematical model of mathematical modeling. We are claiming, then, that category theory mathematically foregrounds, and formalizes, certain observable aspects of the subject of mathematical modeling. 

To make this claim, we must answer the question, \emph{What observable aspect of mathematical models does category theory (CT) foreground}? The answer is, roughly, that CT foregrounds a sense in which each mathematical model is known by its relationship with other mathematical models. That is, CT foregrounds the third postulate as an observable aspect of modeling. It formalizes this postulate in terms of \emph{morphisms}.

Clarifying the above statements will be the subject of the present paper. I will explain how notions in pure mathematics, such as vector spaces or groups, can be viewed as mathematical models, say of linearity or symmetry. I will show how models of linearity, symmetry, and action are all known by the interactions that exist between them. In other words, I am committing to the above description of, and postulates about, models---exemplified by sonar in submarines---and I will present several canonized mathematical concepts, e.g., vector spaces and groups, as mathematical models.

While many mathematicians would agree with the statement that category theory is valuable for understanding and working with mathematical subjects, this paper does not attempt to prove it. However, our second postulate characterizes what such a statement would mean: the value of category theory should be measured by the extent to which it successfully mediates our interactions with mathematical models. If it is valuable, this would imply that category theory foregrounds and faithfully formalizes a propitious aspect of modeling: namely, that to gainfully use models, it is useful to be able to translate between them.

\subsection{Plan of the paper}

Our first goal will be to gently introduce categories in Section~\ref{sec:matrices} using what we hope is a familiar mathematical subject, that of matrix arithmetic. We will then discuss a model of linearity, or flat spaces, in Section~\ref{sec:modeling linearity}, where we will emphasize the categorical perspective, i.e., how the mathematical model of flat spaces is reflected in (and determined by) the rules defining relationships between flat spaces. In Section~\ref{sec:categories}, we define the sort of relationship between categories that captures their structure, namely the functorial relationship. This enables us to consider symmetry and action in Section~\ref{sec:symmetry and action}. Finally, we give a few concluding remarks in Section~\ref{sec:conclusion}.

Throughout the paper we continually return to our three postulates about modeling. In this way, we will be able to view category theory as a mathematical model of mathematical modeling.

\section{Matrices: from groups to enriched categories}\label{sec:matrices}

Our goal in this section is to introduce category theory by considering the case of matrix arithmetic. We will see that all the usual issues regarding dimension and invertibility are actually information about the structure of a category hidden behind the scenes. 

Before we begin, it should be noted that matrices are among the most important tools in mathematical modeling. For example MATLAB, a highly popular technical computing program used by engineers of all kinds, is based primarily on matrix arithmetic. Thus, considering matrices is certainly fair game for thinking about mathematical modeling in the usual sense. Note that the problems usually considered in creating these tools, e.g., speed or accuracy issues, are not being addressed here. Instead, we are considering the abstract idea of matrices. In Section~\ref{sec:matrices review} we will foreground some observable aspects of matrix arithmetic, and in Sections~\ref{sec:monoid matrix mult} -- \ref{subsec:group_enriched} we will formalize them category-theoretically. The result will be a category $\Mat{-}{-}$ that models the subject of matrix arithmetic.

\subsection{Matrices: a review of relevant aspects}\label{sec:matrices review}

A student who is new to linear algebra must learn a few things regarding when matrices can be added and multiplied, the properties of additive and multiplicative identity matrices, the issue of invertibility and non-invertibility, and so on. We now review these because they are precisely what is encoded in the single fact that \emph{matrices form a group-enriched category}. That is, we will introduce category theory by explaining how a certain category models the subject of matrices.

For any natural numbers $m,n\in\NN$, let $\Mat{m}{n}$ denote the set of $m\times n$ matrices---as usual, $m$ is the number of rows and $n$ is the number of columns---and let each entry be a real number. For a matrix $M\in\Mat{m}{n}$, we refer to $(m,n)$ as the \emph{dimension} of $M$ and denote this fact by $\dim(M)=(m,n)$. It is well-known that to add two matrices, say $M+P$, they must have the same dimension. However, to multiply two matrices, say $MP$, there is a different kind of restriction: the middle numbers must agree. More precisely, if $\dim(M)=(m,n)$ and $\dim(P)=(p,q)$, we require $n=p$ in order for the product $MP$ to make sense. In this case, $\dim(MP)=(m,q)$. 

There is a certain matrix $I_n\in\Mat{n}{n}$, called the \emph{$n\times n$ identity matrix}, which looks like this:
\begin{align}\label{dia:identity}
I_n=\begin{bmatrix}
1&0&0&\cdots&0\\
0&1&0&\cdots&0\\
0&0&1&\cdots&0\\
\vdots&\vdots&\vdots&\ddots&\vdots\\
0&0&0&\cdots&1
\end{bmatrix}
\end{align}
A matrix $M$ is called \emph{invertible} if there exists some matrix $N$ such that $MN=I_n$ and $NM=I_n$. Not every matrix is invertible; for example a matrix of all $0$'s, as seen below on the left, is not invertible, but neither is the more average-looking matrix on the right:
$$
\begin{bmatrix}
0&0&0\\
0&0&0\\
0&0&0
\end{bmatrix}
\hspace{.6in}
\begin{bmatrix}
1&2&0\\
0&2&-2\\
-2&-3&-1
\end{bmatrix}$$
Note that for any $m\times n$ matrix $Q$, we have $QI_n=Q$, and for any $n\times p$ matrix $R$ we have $I_nR=R$. 

\subsection{The group of invertible $n\times n$ matrix multiplication}\label{sec:group of invertible matrices}

Many scientists and engineers, from physicists modeling the dynamics of elementary particles, to 3D-animators modeling the changes of camera-angle perspectives on a physical scene, use (either explicitly or implicitly) a \emph{group} of invertible $n\times n$ matrices as part of their toolset. That is, group theory is used in mathematical modeling (in the usual sense). We will see below that group theory is a special case of category theory, and we will explain why groups are showing up in the theory of matrix multiplication. First, however, let us recall what a group is, using matrices as the working example.

Let $\MG{n}$ denote the set of all invertible $n\times n$ matrices. Here are the rules that hold in $\MG{n}$, which make it a \emph{group} in the sense of abstract algebra.
\begin{enumerate}\label{page:defining group}
\item There is an established multiplication formula for $\MG{n}$. In other words, every two elements $M,N\in \MG{n}$ can be multiplied, and the result is again in the group, i.e., $MN\in \MG{n}$. 
\item Multiplying is associative: for any $M,N,P\in \MG{n}$, we have $(MN)P=M(NP)$.
\item There is an established identity element in $\MG{n}$. In other words, there is an element $I_n\in \MG{n}$, such that $I_nM=M=MI_n$ for every $M\in \MG{n}$.
\item There is an established inverse operation in $\MG{n}$. In other words, for every element $M$, there is an established element $N$, often denoted $N=M^{-1}$, such that $MN=I_n=NM$.
\end{enumerate}
These rules encode a notion of symmetry that we will return to later in Section~\ref{sec:modeling with groups}. 

%By symmetry, we only mean an action that can be undone. Each $n\times n$ matrix $M$ encodes an action on the set of $n\times q$ matrix $V$, for any fixed integer $q$. That is, we can multiply $V$ on the left by $M$, and the result will be another $n\times q$ matrix, $MV$. If, as above, we assume that $M$ is invertible, then the action can be undone by multiplying on the left again by the inverse of $M$. This is shown by the following argument:
%$$M^\m1 (MV)=(M^\m1 M)V=I_nV=V.$$
%A \emph{vector} is an $n\times 1$ matrix---a list of $n$ real numbers---which of course is a special case of the above, with $q=1$. In this way, invertible $n\times n$ matrices act as symmetries on vector spaces. We will discuss this more in Section~\ref{sec:modeling linearity}.

\subsection{The monoid of $n\times n$ matrix multiplication}\label{sec:monoid matrix mult}

Recall that $\MS{n}$ denotes the set of all $n\times n$ matrices, including but not limited to the invertible ones. Rather than being a group, $\MS{n}$ is a \emph{monoid}. As such, three out of the four rules
for groups, as enumerated above, are true of $\MS{n}$. Namely,
\begin{enumerate}\label{page:defining monoid}
\item There is an established multiplication formula for $\MS{n}$. In other words, every two elements $M,N\in \MS{n}$ can be multiplied, and the result is again in the monoid, i.e., $MN\in \MS{n}$. 
\item Multiplying is associative: for any $M,N,P\in \MS{n}$, we have $(MN)P=M(NP)$.
\item There is an established identity element in $\MS{n}$. In other words, there is an element $I_n\in \MS{n}$, such that $I_nM=M=MI_n$ for every $M\in \MS{n}$.
\end{enumerate}
A monoid is like a group---elements can be multiplied, multiplication is associative, and there is an identity element---but there is no need for every element of a monoid to be invertible. Not every matrix is invertible, so if we want to think about $n\times n$ matrix multiplication in full generality, we need to use monoids.

A group is thus a special kind of monoid, one in which every element is invertible. In the same way, a monoid is a special kind of \emph{category}, one in which every two elements can be multiplied. Just as we broadened our view from the set of invertible matrices (as a group) to the set of all $n\times n$ matrices (as a monoid), it is now time to further broaden our view to consider all $m\times n$ matrices, i.e., matrices that are not necessarily square. For this we will need a category.

\subsection{The category of matrix multiplication}\label{sec:cat of matrix mult}

Let $\NN$ denote the set of natural numbers. The set of all matrices forms neither a group nor a monoid, but an $\NN$-category, which we will denote $\Mat{-}{-}$. Read the symbol $\Mat{-}{-}$ as ``blank-by-blank matrices." In terms of sets, it is the union
$$\Mat{-}{-}:=\bigcup_{m,n\in\NN}\Mat{m}{n}.$$
The three rules for $\NN$-categories are like those for monoids, but with a slight relaxation in the multiplication rule:%
\footnote{
By $\NN$-category, I mean a category with a position, or \emph{object}, for every natural number $n\in\NN$. For the definition of a general category, first add the following rule before the listed three:
\begin{enumerate}\setcounter{enumi}{-1}
\item There is an established set $\Ob$, whose elements are called \emph{objects}.
\end{enumerate}
Then, replace every occurrence of $\NN$ with an occurrence of $\Ob$. In other words, an $\NN$-category is a category in which $\Ob=\NN$. In the standard definition of a category, $\Mat{m}{n}$ is playing the role of the set of \emph{morphisms} $m\to n$.
}%
\begin{enumerate}\label{page:defining category}
\item There is an established multiplication formula for $\Mat{-}{-}$, which is defined as long as the middle terms agree. In other words, matrices $M\in\Mat{m}{n}$ and $N\in\Mat{p}{q}$ can be multiplied if and only if $n=p$, and the result is again in the category, i.e., $MN\in\Mat{m}{q}$. 
\item Multiplying is associative: for any $M,N,P\in\Mat{-}{-}$, if $MN$ and $NP$ can be multiplied then $(MN)P=M(NP)$. 
\item For each $n\in\NN$ there is an established identity element in $\Mat{n}{n}$. In other words, there is an element $I_n\in\Mat{n}{n}$, such that \mbox{$MI_n=M$} for every $M\in\Mat{m}{n}$ and $I_nN=N$ for every $N\in\Mat{n}{q}$.
\end{enumerate}

Let's consider each dimension $n\in\NN$ to be a kind of context. Then groups are about actions which do not change context and which are reversible; monoids are about actions which do not change context but which may be irreversible; and categories are about actions which may change context and which may be irreversible. While groups and monoids are said to have \emph{elements}, the elements in a category (the elements of $\Mat{-}{-}$ in the above case) are usually called {\em morphisms}. 

Note that although our definition of category looks very much tuned to matrices, it is actually quite general. When someone speaks of a category, they mean nothing more than an establishment of the structures and rules shown in (0)--(3) above.

\subsection{The group-enriched category of matrix arithmetic}\label{subsec:group_enriched}

We still have not grappled with the fact that matrices can be added. For every pair of natural numbers $m,n\in\NN$ the set $\Mat{m}{n}$ can be given the structure of a group, which encodes the addition of $m\times n$ matrices. It is a different group than the one discussed in Section~\ref{sec:group of invertible matrices}---i.e., it encodes addition rather than multiplication---but it is a group nonetheless because it satisfies the same formal rules. That is:
\begin{enumerate}
\item There is an established addition formula for $\Mat{m}{n}$. In other words, every two elements $M,N\in \Mat{m}{n}$ can be added, and the result is again in the group, i.e., $M+N\in \Mat{m}{n}$. 
\item Adding is associative: for any $M,N,P\in \Mat{m}{n}$, we have $(M+N)+P=M+(N+P)$.
\item There is an established (additive) identity element in $\Mat{m}{n}$. In other words, there is an element $Z_{m,n}\in \Mat{m}{n}$, such that $Z_{m,n}+M=M=M+Z_{m,n}$ for every $M\in \Mat{m}{n}$.
\item There is an established (additive) inverse operation in $\Mat{m}{n}$. In other words, for every element $M$, there is an established element $N$, often denoted $N=-M$, such that $M+N=Z_n=N+M$.
\end{enumerate}
Of course, $Z_{m,n}$ is the $m\times n$ matrix of zeros, and $-M$ is the matrix obtained by multiplying each entry in $M$ by $-1$. We saw in Section~\ref{sec:cat of matrix mult} that $\Mat{-}{-}$ is a category; once we include the additive group structure on $\Mat{m}{n}$, our structure becomes a \emph{group-enriched category}.

Thus we see two types of group structures arising in the story of matrix arithmetic: a group encoding matrix addition for $m\times n$ matrices, for any $m,n\in\NN$, and a group encoding matrix multiplication for invertible $n\times n$ matrices, for any $n\in\NN$. And there is further interaction between the additive and multiplicative operations in the category of matrices; namely, multiplication of matrices distributes over addition in the sense that $M(N+P)=MN+MP$. 

The entire addition and multiplication story for matrices, discussed above and in any first course on linear algebra, is subsumed in a single category-theoretic statement: \emph{Matrices form a group-enriched category with objects $\Ob=\NN$.} This articulates:
\begin{itemize}
\item the dimensionality requirements for multiplication and addition,
\item the roles of each identity and zero matrix,
\item the associativity and distributivity laws for multiplication and addition,
\item the existence of additive inverses,
\item the way invertible matrices fit into the picture.
\end{itemize}
Let us clarify the final point. From the perspective that matrices form a category, the notion of invertible matrices comes for free. That is, for every category $\mcC$, and for every object $n$ in it, there is a group of invertible morphisms from $n$ to itself, called the \emph{automorphism group} of $n$. When $\mcC$ is the category $\Mat{-}{-}$, in which each object is a natural number $n\in\NN$, the automorphism group of $n$ is $\MG{n}$. Invertibility is seen as an issue, not about matrices in particular, but about morphisms in any category, and hence the issue is placed in a much broader context.

\section{Modeling linearity}\label{sec:modeling linearity}

In Section~\ref{sec:matrices} we showed how category theory models the subject of matrix arithmetic. But matrices themselves are one aspect of a highly-valued mathematical modeling framework, namely that of vector spaces. Vector spaces are the mathematical model of linearity, as we will discuss in Section~\ref{sec:linearity}. The \emph{category} of vector spaces, $\Vect_\RR$, models this model by foregrounding the sense in which each vector space is known by its relationship with other vector spaces. We will discuss this in Section~\ref{sec:morphisms as structure-keepers}.

\subsection{Vector spaces: the mathematical models of linearity}\label{sec:linearity}

%In his \emph{Discorsi},%
%\footnote{See Heidegger, M. ``Modern Science, Metaphysics, and Mathematics", found in \emph{Basic Writings}.}
%Galileo says, ``I think in my mind of something movable that is left entirely to itself". Here he is doing mathematical modeling---foregrounding, and soon formalizing---the issue of inertia. He elaborates, ``I think of a body thrown on a horizontal plane and every obstacle excluded. This results in what has been given a detailed account in another place, that the motion of the body over this plane would be uniform and perpetual if the plane were extended infinitely."

The notion of linearity shows up in our visual, linguistic, and cognitive interactions with the world. Indeed our visual system is hardwired to highlight straight lines. In our language, simplicity and goodness are often equated with flatness and straightness; English words such as \emph{plain, straightforward, right, direct, correct} and \emph{true} invoke straightness. And linearity also appears to be inherent in our best scientific understanding of the physical universe. For example, general relativity postulates that the universe is locally linearizable (i.e., that close enough to any point, the curved space of the universe can be laid flat), and the predictions founded on that assumption match astoundingly with experiments. Even in mathematics we find that linearity often goes hand in hand with simplicity, where many of our most successful techniques work by reducing a difficult case to a linear one. 

As may by now be clear, when I speak of linearity I am referring not only to lines, but to the general notion of flatness, e.g., to planes and higher-dimensional flat spaces. Unlike in curved spaces, such as soap bubbles, we find that in flat spaces a line can point in any direction without having to curve. In the cases mentioned above, the flat spaces also include a distinguished point, called \emph{the origin}, and the straight line segment from the origin to any other point is called a \emph{vector}. Vectors can be stretched by any scaling factor, and two vectors can be added together; in either case the result is another vector. In mathematics, this kind of abstract flat space is called a \emph{vector space}. 

But what are these vector spaces, and what are their relation with linearity in our visual perception, our language, and our cognition? The relationship here is that vector spaces are valuable \emph{models} of our notions of linearity and flatness. That is,
\begin{enumerate}
\item Vector spaces foreground and formalize certain observed aspects of linearity and relationships between them. 
\item Our visual, linguistic, and cognitive interactions with linearity are successfully mediated by the vector space model.
\item Vector spaces are known by the relationships between them.
\end{enumerate}

Let's begin by considering the observable aspects of linearity. Newton's first law of motion is that the velocity of any object remains constant unless a force is applied to the object. In other words, time acts as a scalar multiplier for the motion of objects: doubling the time doubles the distance traveled but does not alter the direction. Scalar multiplication is at the heart of our notion of linearity: a line is the set of scalar multiples of a given vector. 

But what about higher-dimensional linear spaces? In his thought experiments about motion, Galileo imagined a flat plane on which objects could move unobstructed. A plane can be imagined as a 2-dimensional analogue of a line; it is in some sense the simplest 2-dimensional space. A plane has enough structure to discuss not just distance but also direction. Both the angle between lines in a plane and the degree of inclination of a plane embedded in space were necessary for the laws of motion Galileo wished to discuss. 

There are other aspects to planes and flat spaces that are important in modeling. Namely, in a flat space, the different directions do not interact. That is, moving forward in $x$ does not cause any change to occur in $y$ (compare with a parabola or sphere). An $n$-dimensional vector space is a space in which there are $n$ degrees of freedom, which do not interact with one another. That is, there is a well-defined notion of \emph{coordinate system}, whereby every point in the space is uniquely determined by $n$ numbers. This does not hold on a sphere: while it is the case that every point is determined by its latitude and longitude $(lat,long)$, this determination is not unique. For example, we have $(90^\circ,long_1)=(90^\circ,long_2)$ for every pair of longitudes $long_1,long_2$. In other words, the coordinates of latitude and longitude are not free from one another; they interact at the north and south poles. This issue may seem unimportant, but the point is that such caveats cannot be rectified on the sphere precisely because it is not flat.

As mentioned above, the established mathematical models of flat spaces are vector spaces. A (real) vector space is a collection of vectors, including a vector of length 0, such that each vector can be scaled by any real number, e.g., doubled or tripled in length, and such that any two vectors can be added together to form a new vector. When made precise, these ideas are sufficient to define a coordinate system, or \emph{basis}, which is a minimal set of vectors that span the whole space. For example, if the basis consists of three vectors, $x,y,z$, then every vector $v$ in the space can be uniquely obtained by adding together scalar multiples of the basis vectors, say $v=4x+3y-1.5z$.

Thus the observable aspects of linear spaces foregrounded by the mathematical model of vector spaces are: 
\begin{itemize}
\item the existence of a zero vector, 
\item scalar multiplication of vectors, and 
\item addition of vectors.
\end{itemize}
These aspects and the relationships between them---e.g. commutativity of addition, distributivity of scalar multiplication, existence of additive inverses, etc.---are precisely the content of the formal definition of vector space.%
\footnote{Whenever we speak of vector spaces, we mean finite-dimensional vector spaces over the field $\RR$ of real numbers.}

Ren\'e Descartes (and simultaneously Pierre de Fermat) developed the notion of axes, whose coordinates specify any point in a plane (or 3D space). So coordinate systems were invented long before the abstract notion of vector spaces was. However, a major contribution of vector spaces is that they formalize the ability to change between coordinate systems, e.g. by scaling or adding axes to form new axes. Since one can change coordinates at will, perhaps one does not need them at all. In fact, the vector space concept formalizes the idea, due to Lagrange, that the flat space exists, and calculations can take place, a priori---without need for choosing coordinates.

Geometric considerations, such as how two intersecting lines span a plane and two intersecting planes (in space) form a line, are all completely captured by the vector space model. Moreover, various unexpected exceptions---such as the case where the two intersecting planes happen to be the same (and hence the intersection is not a line but the plane), or the case in which the two planes inhabit a space of more than three dimensions and hence intersect in a point rather than a line---are exposed in the mathematical model. 

We have shown, then, that the vector space model of linearity successfully mediates our visual, linguistic, and cognitive interactions with flat spaces. We have also discussed the system of relationships (commutativity, associativity, distributivity) between various aspects of linearity (zero, scalars, sums). This justifies our first two postulates in this setting, so it remains to explain how vector spaces are known by the relationships between them. It is the job of category theory, as a model, to foreground this third postulate using morphisms. Thus we aim to show that vector spaces (the models of linearity) are known by the morphisms between vector spaces. 

Section~\ref{sec:morphisms as structure-keepers} will kill two birds with one stone: it will show that the vector space model of linearity satisfies the third postulate, and it will show that this satisfaction of the third postulate is itself foregrounded by the category theoretic model $\Vect_\RR$ of the vector space model of linearity.

\subsection{$\Vect_\RR$: the categorical model of vector spaces}\label{sec:morphisms as structure-keepers}

%We will thus consider the observable aspects of vector spaces (zero, scalars, sums) using the mathematical model of category theory. We will highlight how these aspects of a given vector space can be known by the relationships between it and other vector spaces. Category theory formalizes relationships as morphisms (see Section~\ref{sec:cat of matrix mult}). 

Category theory formalizes relationships as morphisms, and the tried-and-true morphisms between vector spaces are called \emph{linear transformations}, which we now recall. Let $V$ and $W$ be vector spaces. A linear transformation $f\taking V\to W$ is a method for converting vectors in $V$ to vectors in $W$. But $f$ cannot be an arbitrary function: in order to be called a linear transformation, it must preserve the structures we have formalized in our model, namely zero, scalars, and sums. That is, it must satisfy
\begin{itemize}
\item $f(0)=0$,
\item $f(r\cdot v)=r\cdot f(v)$, and
\item $f(v+v')=f(v)+f(v')$
\end{itemize}
for any $v,v'\in V$ and $r\in\RR$. In fact, there is a tight connection between the set of linear transformations $V\to W$ and the set $\Mat{|V|}{|W|}$ of $|V|\times|W|$ matrices, where $|V|, |W|\in\NN$ are the dimensions of $V$ and $W$. Once one chooses a basis for $V$ and $W$, this connection becomes a one-to-one correspondence: each matrix represents a single linear transformation. There is an equivalence between the category $\Vect_\RR$ and the category $\Mat{-}{-}$ described in Section~\ref{sec:matrices review}; we will henceforth elide the difference.

The point so far is that the morphisms in the category $\Vect_\RR$ take seriously the structures (zero, scalars, sums) we have formalized in each individual model of linearity. In general, the morphisms in a category are designed to reflect (by preserving) the structures that define each object. Thus, an object in a category is not only (epistemologically) \emph{known by} its relationships with other objects; it (ontologically) \emph{is} what it is by virtue of its relationships with other objects. In a categorical model, the knowing of an object and the being of an object are essentially identical.

In the case of a vector space $V$, pure category theoretic reasoning only allows one to consider those aspects of $V$ that can be defined in terms of morphisms into or out of it. For example, the most important aspect of a vector space is the notion of vector. Any individual vector $v\in V$ defines a ruled line in $V$; that is, all scalar multiples of $v$ lie on a single line, ruled by tick-marks at every integer multiple of $v$. But there is a vector space that represents ruled line-hood itself, namely the one-dimensional vector space $\RR$ of all scalar multiples of $1$. We can now say how the notion of vector is itself defined in terms of morphisms in $\Vect_\RR$: for any vector $v\in V$ there is a unique morphism $\RR\to V$ for which $v$ is the image of the vector $1\in\RR$. So if one wishes to think categorically about vectors in a vector space, one can think only in terms of morphisms in $\Vect_\RR$, i.e., linear transformations. 

In fact, the four notions of vector, zero-vector, scalar multiplication, and addition of vectors can all be understood in terms of linear transformations, as we now explain. First, as mentioned above, a vector in $V$ is the same thing as a linear transformation $\RR\to V$, so we have interpreted single vectors in terms only of morphisms. Second, there is a unique morphism from the 0-dimensional vector space $\RR^0$ to any vector space $V$. There is also a unique morphism $\RR\to\RR^0$, and the composite $\RR\to\RR_0\to V$ picks out the zero-vector in $V$. Third, scaling a vector $v$ by a real number $r\in\RR$ can be understood using only composition of morphisms. It fits beautifully to realize that a scalar like $r$ can be represented as a morphism, namely as the linear transformation $r\colon\RR\to\RR$ sending $x\mapsto rx$. We already know that $v$ is represented by a morphism $\RR\to V$. The scaled vector, $rv\taking\RR\to V$ is then simply the composite of two morphisms $\RR\To{r}\RR\To{v}V$ in $\Vect_\RR$. Fourth, there is an similarly elegant way to understand addition of vectors in terms of morphisms. 

The upshot is that the linear transformations between vector spaces can account for all the formal structure that makes a vector space a vector space. But what about the informal, cognitive aspects of linearity discussed in Section~\ref{sec:linearity}, namely ideas like lines, line segments, inclined planes, projections, intersections, and coordinate systems? We have seen that lines in $V$ are captured by linear transformations $\RR\to V$, and similarly the inclusion of an inclined plane in 3D space is a linear transformation $\RR^2\to\RR^3$. The projection of 3D space onto a line or plane is given by a linear transformation $\RR^3\to\RR$ or $\RR^3\to\RR^2$. Any given coordinate system on a vector space of dimension $n$ is given by a unique linear isomorphism $\RR^n\to V$. The notion of intersecting various lines and planes are also beautifully and completely described by a kind of interplay between morphisms, known as the \emph{pullback}.

It is something of a miracle that so many of our intuitive ideas about vector spaces are describable simply in terms of the morphisms in the category $\Vect_\RR$. However, to a category-theorist, this is routine. Every well-studied category has this property, because this is precisely the observable aspect of modeling that category theory foregrounds. In other words, in order for a subject to be model-able by category theory, its objects of study must be determined by their morphisms to and from other objects. There is a general theorem by Nobuo Yoneda (called Yoneda's lemma) to this effect: any object $c$ in any category $\mcC$ is essentially determined by the morphisms into (or out of) $c$.

In other words, categories can only model ``relationally-determined" subjects, subjects in which each object of study is ontologically determined by its relationships to the others. In this case, what is surprising is how extensive the reach of category theory is. The fact that category theory is consistently used to describe so many parts of mathematics, from topological spaces to groups to ordered sets to measure spaces, means that all of these subjects are relationally-determined. This justifies our third postulate, at least for these specific cases.

In the next section we define functors, which relate different categories like morphisms relate objects. This will prepare us to discuss symmetry and action in Section~\ref{sec:symmetry and action}. 
%We will also give the example of orders as categories, as a simple example. 

\section{Relationships between high-level models}\label{sec:categories}

Throughout this paper, we have roughly been interpreting models and relationships, in the context of category theory, as follows. Each object in a category is a model. For example, each vector space is a model of linearity: $\RR$ is a model of line-hood, $\RR^2$ is a model of plane-hood, etc. But each entire category is also a model, albeit one of a higher level. For example, $\Vect_\RR$ is our model of linearity itself. Each low-level model (object) is defined by its relationships (morphisms) to other low-level models, where these relationships are formalized in the higher-level model (the category). For example, we showed how the linearity---the vector space-ness, as formalized by zero, scalars, and sums---of each vector space $V\in\Vect_\RR$ is defined in terms of morphisms between $V$ and other vector spaces. 

But if low-level models are known by the relationships between them, it should be that high-level models are as well. This is indeed the case. We now discuss the categorical model of higher-level models, i.e., the category of categories.

\subsection{$\Cat$: the category of categories}

There is a larger-sized category, denoted $\Cat$, which includes every normal-sized category as an object.%
\footnote{I apologize for the vagueness with respect to the size issue, but it is not relevant to our discussion. See Shulman, ``Set theory for category theory" for details.} 
If $\mcC$ and $\mcD$ are categories, a morphism between them in $\Cat$ is called a \emph{functor}, and can be denoted $F\taking\mcC\to\mcD$. 

Just like any morphism, a functor is a relationship between two categories that preserves the defining structure of categories. Recall from Section~\ref{sec:cat of matrix mult} that the defining structure of a category $\mcC$ is a set of objects, a set of morphisms, an identity morphism for each object, and a formula for composing morphisms. These satisfy some laws, namely that composing any morphism $g$ with an identity morphism gives back $g$ and that composition is associative. To say that functors preserve the structures that define categories is shorthand for the following more formal statement. A functor $F\taking\mcC\to\mcD$ must
\begin{itemize}
\item assign to each object $c\in\Ob(\mcC)$ an object $F(c)\in\Ob(\mcD)$,
\item assign to each morphism $g\taking c\to c'$ in $\mcC$ a morphism $F(g)\taking F(c)\to F(c')$ in $\mcD$,
\item ensure that the morphism assigned to each identity in $\mcC$ is an identity in $\mcD$, i.e., for all $c\in\Ob(\mcC)$, a functor $F$ must ensure that $F(\id_c)=\id_{F(c)}$, and
\item ensure that the composition formula in $\mcC$ is compatible with the composition formula in $\mcD$, i.e., $F(g\circ_\mcC h)=F(g)\circ_\mcD F(h)$.
\end{itemize}

There is a very basic category, which represents object-hood in $\Cat$ the way that $\RR$ represents ruled line-hood in $\Vect_\RR$. Let $\bigob\in\Ob(\Cat)$ denote the category with one object, $\Ob(\bigob)=\{\bullet\}$; only one morphism, $\id_\bullet\taking \bullet\to\bullet$; and the composition formula says $\id_\bullet\circ\id_\bullet=\id_\bullet$. This category $\bigob$ might be pictured 
\begin{align}\label{dia:objecthood}
\fbox{$\bullet$}
\end{align}
For any category $\mcC$, the objects of $\mcC$ are the same as the functors $\bigob\to\mcC$. Similarly, there is a category $\bigmor$, which might be pictured
\begin{align}\label{dia:morphismhood}
\fbox{$\bullet\too\bullet$}
\end{align} 
that represents morphism-hood in $\Cat$. That is, for any category $\mcC$, the morphisms of $\mcC$ are the same as functors $\bigmor\to\mcC$. There is also a functor that represents identity morphisms as well as a functor that represents composition formulas. 

In other words, the morphisms in $\Cat$ can account for all the formal structure that makes a category a category: objects, morphisms, identities, and compositions. This justifies our third postulate in the case of categories: the structure of any given category $\mcC$ is completely determined by the system of functors that map to it from other categories.

\subsection{Two justifications for the third postulate in general}\label{sec:third postulate}

Our goal has been to show that category theory is a mathematical model of mathematical modeling. We discussed the rough meaning of this statement by introducing three postulates about modeling in Section~\ref{sec:intro}, and we said that category theory foregrounds and formalizes the third: that each model is known by its relationships with other models. We have also justified this postulate in several cases; we now justify it in general. 

In a single category, Yoneda's lemma offers one explication and formal justification of the third postulate, given the interpretation of models and relationships, laid out in the first paragraph of Section~\ref{sec:categories}. However, any functor between two categories provides another way to relate these high-level models, namely by something reminiscent of analogy. That is, if $F\taking\mcC\to\mcD$ is a functor, it relates each object $c\in\Ob(\mcC)$ with an object $F(c)\in\Ob(\mcD)$. The morphisms between $c$ and its neighbors in $\mcC$ are preserved by $F$, which makes $F$ act as a sort of analogy, as observed by Brown and Porter.%
\footnote{See Brown, R.; Porter, T. \emph{Category theory: an abstract setting for analogy and comparison}, \url{http://pages.bangor.ac.uk/~mas010/pdffiles/Analogy-and-Comparison.pdf}, whose aim is similar to, but whose scope is more ambitious than, the present paper. }
Our third postulate was intended as a ordinary-language hybrid of the Yoneda concept and the functorial analogy concept.

We now return to a few of the most important areas of mathematical modeling in the usual sense: symmetry and dynamics. The Yoneda interpretation of our third postulate is always at work, but we will briefly explore the second interpretation: how high-level analogies bring out the features of a subject.

\section{Symmetry and action}\label{sec:symmetry and action}

Consider the reflective symmetry of the visible human body or the rotational symmetry of Escher's \emph{Drawing Hands}. The question of symmetry is about reversible action-ability, e.g., the ability to reflect an object and the ability to rotate an object by $180^\circ$ are both reversible. To say that an object is symmetrical with respect to an act is to say that it does not change when it undergoes this act; e.g., the visible human body is unchanged as it undergoes mirror reflection and \emph{Drawing Hands} is unchanged as it undergoes $180^\circ$ rotation. 

We begin in Sections~\ref{sec:modeling with groups}~and~\ref{sec:monoids} by considering various questions of symmetry and action-abilities in the above sense, but this is mainly to set up our main point. Namely, a group $G$ of symmetries (reversible action-abilities) exists independently of any thing that is so-symmetrical (unchanged under corresponding actions). But then what is the connection of the symmetry to the symmetrical? It is captured simply by a functorial connection between one category and another, namely the reversible action category $G$ and the space category $\Vect_\RR$. The functorial connection between a category like $G$%
\footnote{By the phrase ``a category like $G$", we mean a category with one object; see Section~\ref{sec:modeling with groups}.} 
and another category $\mcC$ is called an \emph{action} of $G$ on an object of $\mcC$. We will study actions in Section~\ref{sec:action}. Finally in Section~\ref{sec:dynamical}, we will consider dynamical systems, which are classical exemplars of mathematical modeling, and which again are nothing but functorial connections (roughly, between a time category and a space category).

\subsection{Modeling reversible action-ability with groups}\label{sec:modeling with groups}

A square $S$, with unit-length sides, centered at the origin in $\RR^2$ is symmetrical with respect to eight action-abilities: it can be rotated and reflected in any combination:
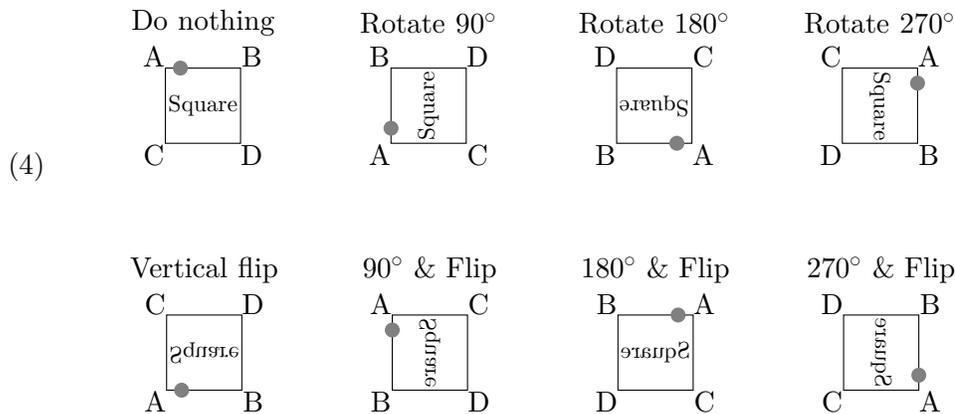
\begin{align}\label{dia:square}
\begin{tikzpicture}
%Do nothing
	\draw (.5,1.6) node {Do nothing};
	\draw (0,0) rectangle (1,1);
	\draw(.5,.5) node [xscale=.8,yscale=.8]{Square};
	\draw (0.2,1) node [circle,fill=gray,inner sep=2pt] {};
	\draw (-.15,-.15) node {C};
	\draw (-.15,1.15) node {A};
	\draw (1.15,-.15) node {D};
	\draw (1.15,1.15) node {B};
%Rotate $90^\circ$
	\draw (3.5,1.6) node {Rotate $90^\circ$};
	\draw (3,0) rectangle (3+1,1);
	\draw(3.5,.5) node [rotate=90,xscale=.8,yscale=.8]{Square};
	\draw (3,.2) node [circle,fill=gray,inner sep=2pt] {};
	\draw (3+-.15,-.15) node {A};
	\draw (3+-.15,1.15) node {B};
	\draw (3+1.15,-.15) node {C};
	\draw (3+1.15,1.15) node {D};
%Rotate $180^\circ$
	\draw (6+.5,1.6) node {Rotate $180^\circ$};
	\draw (6+0,0) rectangle (6+1,1);
	\draw(6+.5,.5) node [rotate=180,xscale=.8,yscale=.8]{Square};
	\draw (6.8,0) node [circle,fill=gray,inner sep=2pt] {};
	\draw (6+-.15,-.15) node {B};
	\draw (6+-.15,1.15) node {D};
	\draw (6+1.15,-.15) node {A};
	\draw (6+1.15,1.15) node {C};
%Rotate $270^\circ$
	\draw (9+.5,1.6) node {Rotate $270^\circ$};
	\draw (9+0,0) rectangle (9+1,1);
	\draw(9+.5,.5) node [rotate=270,xscale=.8,yscale=.8]{Square};
	\draw (10,.8) node [circle,fill=gray,inner sep=2pt] {};
	\draw (9+-.15,-.15) node {D};
	\draw (9+-.15,1.15) node {C};
	\draw (9+1.15,-.15) node {B};
	\draw (9+1.15,1.15) node {A};
\end{tikzpicture}
\\\nonumber\\\nonumber
\begin{tikzpicture}
%Vertical flip
	\draw (.5,1.6) node {Vertical flip};
	\draw (0,0) rectangle (1,1);
	\draw(.5,.5) node [xscale=.8,yscale=-.8]{Square};
	\draw (0.2,0) node [circle,fill=gray,inner sep=2pt] {};
	\draw (-.15,-.15) node {A};
	\draw (-.15,1.15) node {C};
	\draw (1.15,-.15) node {B};
	\draw (1.15,1.15) node {D};
%Flip and $90^\circ$
	\draw (3+.5,1.6) node {$90^\circ$ \& Flip};
	\draw (3+0,0) rectangle (3+1,1);
	\draw(3+.5,.5) node [rotate=90,xscale=-.8,yscale=.8]{Square};
	\draw (3,.8) node [circle,fill=gray,inner sep=2pt] {};
	\draw (3+-.15,-.15) node {B};
	\draw (3+-.15,1.15) node {A};
	\draw (3+1.15,-.15) node {D};
	\draw (3+1.15,1.15) node {C};
%Flip and $180^\circ$
	\draw (6+.5,1.6) node {$180^\circ$ \& Flip};
	\draw (6+0,0) rectangle (6+1,1);
	\draw(6+.5,.5) node [rotate=180,xscale=.8,yscale=-.8]{Square};
	\draw (6.8,1) node [circle,fill=gray,inner sep=2pt] {};
	\draw (6+-.15,-.15) node {D};
	\draw (6+-.15,1.15) node {B};
	\draw (6+1.15,-.15) node {C};
	\draw (6+1.15,1.15) node {A};
%Flip and $270^\circ$
	\draw (9+.5,1.6) node {$270^\circ$ \& Flip};
	\draw (9+0,0) rectangle (9+1,1);
	\draw(9+.5,.5) node [rotate=270,xscale=-.8,yscale=.8]{Square};
	\draw (10,.2) node [circle,fill=gray,inner sep=2pt] {};
	\draw (9+-.15,-.15) node {C};
	\draw (9+-.15,1.15) node {D};
	\draw (9+1.15,-.15) node {A};
	\draw (9+1.15,1.15) node {B};
\end{tikzpicture}
\end{align}
Neither the word \emph{Square}, the gray dot, nor the labels $A, B, C, D$ are symmetrical with respect to the same eight-element group that the square itself is. They are drawn in (\ref{dia:square}) to show that unlike the square, they undergo change when we act on them in these eight ways.

Every one of these eight actions on the square is reversible, and its reverse is another one of the eight actions. The actions are also serializable, in the sense that if each of $a_1,\ldots,a_n$ is one of the eight, then so is the process obtained by doing them in series, denoted $a_1a_2\cdots a_n$. In Section~\ref{sec:group of invertible matrices} we called this serialization \emph{multiplication}, but in category-theoretic terms it is called \emph{composition}. The above eight-element set, which has an identity element and is closed under inverses and compositions, is a group (see Section~\ref{sec:group of invertible matrices}), called the \emph{dihedral group of order 8}, and denoted $D_8$. 

Note that the group $D_8$ acts on the square, but in fact $D_8$ exists independently of the square. That is,  $D_8$ could just as well act on an octagon or on a single point at the origin. We refer to the elements of $D_8$ as \emph{action-abilities} because each is able to act on a variety of things. The elements only become \emph{actions} when they are applied to something, such as a square. We will discuss the notion of action itself in Section~\ref{sec:action}.

Every group can be modeled as a category $\mcG$ with many morphisms, each corresponding to an action-ability, but with only one object. The unique object of $\mcG$ stands for ``the abstract thing that is unchanged under these actions." The identity morphism corresponds to the ability not to act, and the composition formula in $\mcG$ corresponds to the requirement that the serialization of action-abilities is an action-ability. Groups are categories that encode symmetries, which we define as action-abilities that can be undone. That is, for every morphism $a\in\mcG$, there is a morphism $b$ for which the serializations $a$-then-$b$ and $b$-then-$a$ are equal to the identity (non-)action.

The group $D_8$ is one type of symmetry; there are many others. For example, consider the line that goes through the origin in $\RR^2$ at an angle of $30^\circ$. If you take each point on the line and multiply its distance from the origin by a non-zero number $k$, but do not change the angle, the total result is again the same $30^\circ$ line: it is unchanged by non-zero scaling. Thus there is a different group, $\RR_{\neq 0}$, of non-zero scaling abilities. It is a different type of symmetry than $D_8$, but just as much able to encode a kind of ability to act reversibly.

Before we move on, note that we now know that each group is a category. A functor between two groups is called a \emph{group homomorphism}. The category of all groups and group homomorphisms, denoted $\Grp$, is the high-level model of symmetry itself, just like the category $\Vect_\RR$ is the high-level model of linearity.

\subsection{Modeling action-ability with monoids}\label{sec:monoids}

As mentioned in Section~\ref{sec:monoid matrix mult}, a monoid is like a group, except that not all its elements need be invertible. For example, there is a monoid $M$ consisting of four action-abilities, which we can depict using their action on a windowpane. 
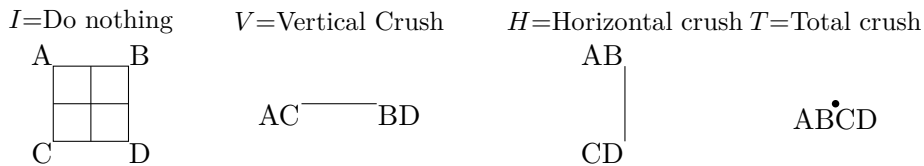
\begin{align*}
\begin{tikzpicture}
%Do nothing
	\draw (.5,1.6) node {\small$I$=Do nothing};
	\draw (0,0) rectangle (1,1);
	\draw (.5,0) -- (.5,1);
	\draw (0,.5) -- (1,.5);
	\draw (-.15,-.15) node {C};
	\draw (-.15,1.15) node {A};
	\draw (1.15,-.15) node {D};
	\draw (1.15,1.15) node {B};
%Vertical crush
	\draw (3.3+.5,1.6) node {\small$V$=Vertical Crush};
	\draw (3.3+0,0.5) -- (3.3+1,0.5);
	\draw (3.3+-.3,-.15+.5) node {AC};
	\draw (3.3+1.3,-.15+.5) node {BD};
%Horizontal crush
	\draw (6.6+.5+.5,1.6) node {\small$H$=Horizontal crush};
	\draw (6.6+.5+0+.5,0) -- (6.6+.5+0+.5,1);
	\draw (6.6+.5+-.3+.5,1.15) node {AB};
	\draw (6.6+.5+-.3+.5,-.15) node {CD};
%Total crush
	\draw (9.9+.5,1.6) node {\small$T$=Total crush};
	\draw (9.9+.5,0+.5) node [circle,fill=black,inner sep=1pt] {};
	\draw (9.9+.5,.3) node {ABCD};
\end{tikzpicture}
\end{align*}
Any sequence of these action-abilities can be serialized, and the result is also one of these four action-abilities. Also, there is an identity, ``do nothing" action-ability. Thus we have a monoid, and it is not a group because three of the action-abilities are irreversible. Each element of $M$ acts on the windowpane: it sends every point in the windowpane to another point in the windowpane.

A monoid is a category with one object, but with possibly many morphisms from that object to itself, as depicted here: 
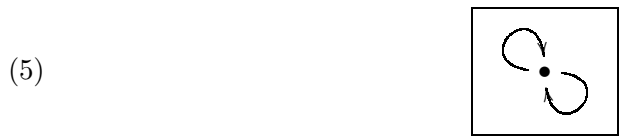
\begin{align}\label{dia:monoid and order}
\parbox{.75in}{
\fbox{\xymatrix{
\bullet\ar@(l,u)[]^~\ar@(r,d)[]^~
}
}}
\end{align}
The fact that there is only one object means that any two morphisms can be composed; this is the serializability. 

In passing we note that the connection between groups and monoids (every group is a monoid but not vice versa) is captured by a high-level analogy, i.e. a functor $\Grp\to\Mon$, where $\Mon$ is the category of all monoids. There are even functors going the other way, $\Mon\to\Grp$, but we will not discuss any of them here.
%In the above case, the composition of any element $x$ with the identity $I$ or with $x$ itself  returns $x$. The composition of any other two elements returns the ``Total crush" element. For example, vertically crushing the windowpane and then horizontally crushing the result is the same as totally crushing it. Here is the composition formula in the form of a multiplication table:
%$$\small
%\begin{array}{l | l l l l}
%\circ&I&V&H&T\\\hline
%I&I&V&H&T\\
%V&V&V&T&T\\
%H&H&T&H&T\\
%T&T&T&T&T
%\end{array}
%$$

\subsection{Modeling action with outgoing functors}\label{sec:action}

The four elements of the monoid $M$ from the previous section (\ref{sec:monoids}) can be represented by the following matrices:
\begin{align}\label{dia:representation}
I\mapsto\begin{bmatrix}1&0\\0&1\end{bmatrix}\hsp
V\mapsto\begin{bmatrix}1&0\\0&0\end{bmatrix}\hsp
H\mapsto\begin{bmatrix}0&0\\0&1\end{bmatrix}\hsp
T\mapsto\begin{bmatrix}0&0\\0&0\end{bmatrix}
\end{align}
But how exactly can we enunciate this connection between elements of our monoid $M$ and these matrices, i.e., morphisms in $\Vect_\RR$.

The most concise and straightforward way I know is to use functors. A functor $F\taking M\to\Vect_\RR$ assigns to each object in $M$ an object in $\Vect_\RR$. Since there is only one object in $M$, we get only one vector space; in the above case, it is $\RR^2$. A functor $F$ also assigns to each morphism in $M$ a morphism in $\Vect_\RR$; the four elements $\{I,V,H,T\}$ of the monoid are then sent to four linear transformations $\RR^2\to\RR^2$. These are represented by the four matrices above in \eqref{dia:representation}. Moreover, the composition formula for $M$ is preserved by $F$. This ensures that if an equation holds in $M$, e.g., $VH=T$, then the corresponding equation holds in $\Vect_\RR$, e.g., 
$$
\begin{bmatrix}1&0\\0&0\end{bmatrix}\begin{bmatrix}0&0\\0&1\end{bmatrix}=\begin{bmatrix}0&0\\0&0\end{bmatrix}
$$

In mathematics, if $M$ is a monoid (or a group), we define an \emph{$M$-action on a set} to be a functor $M\to\Set$; we define an \emph{$M$-action on a vector space} to be a functor $M\to\Vect_\RR$; and so on. Thus in (\ref{dia:representation}) we have established an $M$-action on $\RR^2$. Similarly, the group $\MG{n}$ of invertible $n\times n$ matrices, discussed in Section~\ref{sec:group of invertible matrices} acts on $\RR^n$. The study of functors from a group to $\Vect_\RR$ is often called \emph{representation theory} in mathematics. 

There is a natural category structure on the class of functors between any two categories. These functor categories are often important, as they are in representation theory. In the next section we briefly describe dynamical systems in these terms.

\subsection{Dynamical systems}\label{sec:dynamical}

To model the behavior of a system that changes in time, one must decide whether to formalize its change as continuous or discrete. For example, the internal states of a computer may be best modeled with a discrete dynamical system, if it has an internal clock that dictates discrete times at which changes can occur. On the other hand, the concentrations of chemicals in a reaction process are changing continuously, so the reaction is often modeled by a continuous dynamical system. Here we will focus on autonomous dynamical systems.

A \emph{discrete dynamical system} is defined to be a set $S$ and a \emph{where-to-go-next} function $f\taking S\to S$. A case where $S$ has eleven elements may be depicted
$$
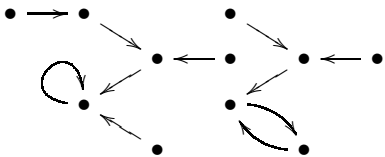
\xymatrix@R=5pt@C=15pt{
\bullet\ar[r]&\bullet\ar[dr]&&\bullet\ar[dr]&\\
&&\bullet\ar[dl]&\bullet\ar[l]&\bullet\ar[dl]&\bullet\ar[l]\\
&\bullet\ar@(l,u)[]&&\bullet\ar@/^.5pc/[dr]\\
&&\bullet\ar[ul]&&\bullet\ar@/^.5pc/[ul]
}
$$

In fact a discrete dynamical system is the same thing as a functor $\NN\to\Set$, where $\NN$ is the monoid of natural numbers under addition. In other words, consider $\NN$ as a one-object category with morphisms $0,1,2,\ldots$, where  composing $i$ with $j$ is given by the formula $i+j$. A functor $\NN\to\Set$ consists of a set $S$ and a function $f^n\taking S\to S$ for every natural number $n$. The fact that a functor must preserve the composition formula implies that $f^i(f^j(s))=f^{i+j}(s)$ for every $s\in S$.

A \emph{continuous dynamical system} (sometimes called a \emph{flow}) is a topological space $S$ and a continuous function $f\taking S\times\RR_{\geq 0}\to S$, such that 
$$f(s,0)=s\hsp\tn{and}\hsp f(f(s,t_1),t_2)=f(s,t_1+t_2).$$
Here $f(s,3.14)$ would tell us where the point $s$ will be after 3.14 units of time. In fact, a continuous dynamical system can be modeled categorically as a topologically-enriched functor $\RR_{\geq 0}\to\Top$. The reader is not expected to understand this statement exactly, but the idea is that we can concisely capture the definition of dynamical systems using functors. 

We have now shown that the class of dynamical systems (either discrete or continuous) is itself a \emph{category of functors} from a time category to a space category. For discrete dynamical systems, time and space are modeled discretely (with time as $\NN$ and space as a set). For continuous dynamical systems, they are modeled continuously (with time as $\RR_{\geq 0}$ and space as a topological space). However, the dynamical system itself is a functor between these categories.

Time and space can be modeled as independent categories, but part of the human concept of time is that it acts on objects in space. That is, we understand our model of time and our model of space by connecting the two. This is an example of the analogical interpretation of the third postulate, and we have modeled it formally in this section using functors.

\section{Conclusion}\label{sec:conclusion}

Karl Popper said, ``A theory that explains everything, explains nothing." If category theory models algebra, geometry, logic, computer science, probability, and more, is it not trying to be a model of everything? The point is a bit subtle.

Category theory is not a theory of everything. It is more like, as topologist Jack Morava put it,%
\footnote{See Morava, J. ``Theories of anything", \url{http://arxiv.org/abs/1202.0684}.}
``a theory of theories of anything". In other words, it is a model of models. It leaves each subject alone to solve its own problems, to sharpen and refine its toolset in the ways it sees fit. That is, CT does not micromanage in the affairs of any discipline. However, describing any discipline categorically tends to bring increased conceptual clarity, because conceptual clarity is CT's main concern, its domain of expertise. And it does enforce certain principles; for example the theorem of Yoneda, discussed above, ensures that it is not the objects, but rather the relationships between objects, that determine the essence of any category. Finally, category theory allows one to compare different models, thus carrying knowledge from one domain to another, as long as one can construct the appropriate ``analogy", i.e., functor.

Category theory has continually sharpened and refined its own toolset, i.e., its ability to articulate the various objects, relationships, properties, structures, and methods that show up throughout mathematics. It consistently shows itself as a powerful mode of mathematical thinking, and there is no a priori reason it cannot be similarly successful in science more broadly. 

However, a question remains: what of the modeler? Who is the one that decides that a certain category adequately models a certain subject? Who is the one that finds value in the category-theoretic mode of thought? Perhaps category theory can aid in a mathematically rigorous form of phenomenological reduction,%
\footnote{See Cogan, J., ``The phenomenological reduction", \url{http://www.iep.utm.edu/phen-red/}.}
in which the process of thinking is itself elucidated. Category theory could be considered a truly profound model of modeling if it could model the cognitive apparatus itself, i.e., if it could foreground and mathematically communicate the relationship between subject, model, and modeler.

\subsection*{Acknowledgements}
I want to thank Allen Brown, Patrick Schultz, Tsung-Yun Tzeng, and Dmitry Vagner, as well as the referees, for their helpful comments on various drafts of this paper. This project was supported by ONR grant N000141310260 and AFOSR grant FA9550-14-1-0031.

\end{document}